\documentclass{article}
\usepackage{amssymb,amsmath}
\usepackage{graphics}

\def\qed{\hfill {\hbox{${\vcenter{\vbox{               
   \hrule height 0.4pt\hbox{\vrule width 0.4pt height 6pt
   \kern5pt\vrule width 0.4pt}\hrule height 0.4pt}}}$}}}

\def\tr{\triangleright}

\newtheorem{proposition}{Proposition}[section]

\newtheorem{corollary}{Corollary}[section]
\newtheorem{example}{Example}[section]
\newtheorem{conjecture}{Conjecture}[section]

\newenvironment{proof}[1][Proof]{\medskip\noindent{\bf #1.}\quad}%
{\qed\par\medskip}

\author{
Xiao-Song Lin\footnote{Partially supported by NSF.} \\
\texttt{\small xl@math.ucr.edu} 
\and
Sam Nelson \\
\small{\texttt{knots@esotericka.org}}
\and
\begin{tabular}{c}
\small{Department of Mathematics, University of California, Riverside} \\
\small{900 University Avenue, Riverside, CA, 92521}
\end{tabular}
}

\title{\Large \textbf{On Generalized Knot Groups}} 

\date{} 

\begin{document}

\maketitle

\abstract{Generalized knot groups $G_n(K)$ were introduced first by 
Wada and Kelly independently. The classical knot group is the first one 
$G_1(K)$ in this series of finitely presented groups. For each natural 
number $n$, $G_1(K)$ is a subgroup of $G_n(K)$ so the generalized knot 
groups can be thought of as extensions of the classical knot group. For the 
square knot $SK$ and the granny knot $GK$, we have an isomorphism 
$G_1(SK)\cong G_1(GK)$. From the presentations of $G_n(SK)$ and $G_n(GK)$, 
for $n>1$, it seems unlikely that $G_n(SK)$ and $G_n(GK)$ would be 
isomorphic to each other. We are able to show that for many finite 
groups $H$, the numbers of homomorphisms from $G_n(SK)$ and $G_n(GK)$ to $H$, 
respectively, are the same. Moreover, the numbers of conjugacy classes of 
homomorphisms from $G_n(SK)$ and $G_n(GK)$ to $H$, respectively, are also 
the same. It remains a challenge to us to show, as we would like to 
conjecture, that $G_n(SK)$ and $G_n(GK)$ are not isomorphic to each other 
for all $n>1$.  }

\section{\large \textbf{Introduction}}

Generalized knot groups were introduced by Wada \cite{W} and 
Kelly \cite{K} independently. Wada arrived at these group invariants of knots
by searching for homomorphisms of the braid group $B_n$ 
into $\text{Aut}(F_n)$, while Kelley's work was related to knot quandles and 
Wirtinger-type presentations.  

We find it convenient to introduce generalized knot groups through 
the language of quandles (see \cite{FR} and \cite{J}). A \textit{quandle} is 
a set $Q$ with a binary operation $\tr:Q\times Q \to Q$ satisfying the 
conditions
\newcounter{q}
\begin{list}{(\roman{q})}{\usecounter{q}}
\item for all $x\in Q$, $x\tr x = x$,
\item for all $x,y\in Q$, there is a unique $z\in Q$ with $x=z\tr y$, 
and
\item for all $x,y,z\in Q$, $(x\tr y)\tr z= (x\tr z)\tr (y\tr z)$.
\end{list}
If $(Q,\tr)$ satisfies $(ii)$ and $(iii)$ but not necessarily $(i)$, then 
$(Q,\tr)$ is a \textit{rack}.

Quandles and knots are closely related; if we interpret quandle elements as 
arcs in a link diagram, then the three quandle axioms are just the three 
Reidemeister moves. Quandles and groups are also closely related; indeed, a 
group is a quandle with quandle operation given by conjugation. 

More precisely, given a group $G$, there is a quandle $\mathrm{Conj}(G)$ 
with underlying set $G$ and quandle operation $\tr:G\times G\to G$ given by
$x\tr y= y^{-1}xy$. Conversely, for any quandle $Q$ there is a group called
the \textit{associated group} of $Q$, 
$\mathrm{As}(Q)=F(Q)/\langle y^{-1}xy(x\tr y)^{-1} \ \forall x,y\in Q \rangle$,
that is, the free group on $Q$ modulo the normal subgroup generated by 
relations obtained by setting $y^{-1}xy$ equal to $x\tr y$ for all elements 
of $Q$. $\mathrm{As}(Q)$ is also called $\mathrm{Adconj}(Q)$, since the functor
$\mathrm{As}:\mathbf{QUANDLES}\to\mathbf{GROUPS}$ is the left adjoint to
the functor $\mathrm{Conj}:\mathbf{QUANDLES}\to\mathbf{GROUPS}.$

It is easy to check that for any $n\in \mathbb{Z}$, 
$\mathrm{Conj_n}:\mathbf{GROUPS}\to\mathbf{QUANDLES}$
defined by $\mathrm{Conj_n}(G)=(G,x\tr y= y^{-n}xy^n)$ is likewise a 
functor from the category of groups to quandles; its left adjoint, 
$\mathrm{As_n}:\mathbf{QUANDLES}\to\mathbf{GROUPS}$ given by 
$\mathrm{As_n}(Q)=F(Q)/\langle x^{-n}yx^n(x\tr y)^{-1} \ 
\forall x,y\in Q \rangle$ is then a functor from quandles to groups. In 
particular, the fundamental group of a link complement is the associated 
group of the fundamental quandle of the link complement; if we now consider 
the $n$th associated group of the link's fundamental quandle, we obtain a 
new group invariant of links, which we will denote by $G_n(K)$ for a link 
$K$. Since $G_1(K)$ is just the usual fundamental group of $S^3\setminus K$, 
these groups form a family of generalized knot groups.

It is well known that the knot group $G_1(K)$ alone is not strong enough to 
classify knots up to ambient isotopy. So, a natural question is whether the 
isomorphism types of generalized knot groups could be used as classifying 
invariants of knots. Unfortunately, this is not the case, as proposition 
\ref{trefoil} shows. However, we suspect that generalized knot groups hold
additional information about knot type that is not present in the
usual fundamental group.

As a testing case, we consider the square knot $SK$ and the granny knot
$GK$. Then $G_1(SK)$ and $G_1(GK)$ are isomorphic to each other. To check 
whether $G_n(SK)$ and $G_n(GK)$ are isomorphic to each other, for $n\geq 2$,
we programmed our computer to calculate the numbers of homomorphisms of 
$G_n(SK)$ and $G_n(GK)$ into a finite group $H$, respectively. To our 
surprise, for the many finite groups that we tested, these numbers are 
always equal to each other. We will show that this is indeed the case for 
any finite group satisfying a certain property. Furthermore, we will show 
that the numbers of conjugacy 
classes of homomorphisms of $G_n(SK)$ and $G_n(GK)$ into $H$, respectively, 
are also the same. Thus, we can not distinguish these two groups $G_n(SK)$ 
and $G_n(GK)$, for each $n\geq 2$, by simply counting (conjugacy classes
of) homomorphisms into any finite group satisfying a certain property.  
 
From the presentations of $G_n(SK)$ and $G_n(GK)$, it is very tempting to 
conjecture that $G_n$ distinguishes the square knot from the granny knot. 
If this is the case, then the groups $G_n$ contain additional information 
about knot type not contained in the fundamental group. Unfortunately, the 
situation turned out to be much more subtle than we originally thought. See 
the discussion in the last section about some computationally intensive 
attempts to show that the groups $G_n(SK)$ and $G_n(GK)$ are not isomorphic 
for $n\geq2$. This conjecture remains open.

\section{\large \textbf{The fundamental quandle and $G_n$}}

Let $L$ be a link diagram. The fundamental quandle of $L$ has presentation
\[Q(L)=\langle x_1, \dots ,x_m \ | \ x_i\tr x_j=x_k \ \mathrm{or} 
\ x_k\tr x_j=x_i \ \ \ i=1, \dots, m \rangle\]
with one generator $x_i$ for each arc in $L$ and one relation
for each crossing, either of the form $x\tr y= z$  where $x$ is the incoming 
underarc, $y$ is the overarc and $z$ is the outgoing overarc at a positive 
crossing, or $x\tr y= z$  where $x$ is the outgoing underarc, $y$ is the 
overarc and $z$ is the incoming overarc at a negative crossing.

\begin{figure}[!ht]
\[ \includegraphics{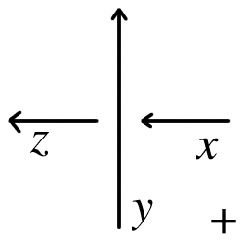} \ \ \ \includegraphics{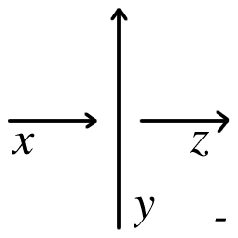} \]
\[  x\tr y=z  \ \ \ \ \ \ \ \ \ \ \ \ \ z\tr y= x\]
\[  x=y^nzy^{-n}\ \ \ \ \ \ \ \ \ \ \ x=y^{-n}zy^{n} \]
\caption{Quandle and $G_n$ relations at a crossing.}
\label{qr}
\end{figure}

This then permits us to give a combinatorial description of $G_n(L)$ defined 
from a link diagram, namely $G_n(L)$ is the group with presentation
\[G_n(L)=\langle x_1, \dots ,x_m \ | \ x_i=x_j^nx_kx_j^{-n} \ \mathrm{or} 
\ x_i=x_j^{-n}x_kx_j^n, \ i=1,\dots, m\rangle, \] with the type of relation 
determined by the sign of the crossing.

\begin{proposition} \label{trefoil}
The right- and left-handed trefoils have isomorphic generalized groups
$G_n$ for all $n\in \mathbb{Z_+}$.
\end{proposition}

\begin{proof}
The right-hand trefoil has three positive crossings and generalized knot group
\[G_n(T_r)=\langle a,b,c \ | \ a=b^ncb^{-n}, b=c^nac^{-n}, c=a^nba^{-n} 
\rangle. \] Using the fact that $(x^nyx^{-n})^m=x^ny^mx^{-n}$, we can 
eliminate the generator $c$ to obtain
\begin{eqnarray*}
G_n(T_r) & = & \langle a,b,| \ a=b^n(a^nba^{-n})b^{-n},
 b=(a^nba^{-n})^na(a^nba^{-n})^{-n} \rangle \\
& = & \langle a,b,| \ a=b^n(a^nba^{-n})b^{-n},
 b=(a^nb^na^{-n})a(a^nb^{-n}a^{-n}) \rangle \\
& = & \langle a,b,| \ a=b^na^nba^{-n}b^{-n}, b=a^nb^nab^{-n}a^{-n} \rangle \\
& = & \langle a,b,| \ ab^na^n=b^na^nb, ba^nb^n=a^nb^na \rangle.
\end{eqnarray*}

On the other hand, the left-handed trefoil has three negative crossings 
and generalized knot group
\[G_n(T_l)=\langle a,b,c \ | \ a=b^{-n}cb^{n}, b=c^{-n}ac^{n}, c=a^{-n}ba^{n} 
\rangle. \]
Thus, 
\begin{eqnarray*}
G_n(T_l) & = & \langle a,b,| \ a=b^{-n}(a^{-n}ba^{n})b^{n},
 b=(a^{-n}ba^{n})^{-n}a(a^{-n}ba^{n})^{n} \rangle \\
& = & \langle a,b,| \ a=b^{-n}(a^{-n}ba^{n})b^{n},
 b=(a^{-n}b^{-n}a^{n})a(a^{-n}b^{n}a^{n}) \rangle \\
& = & \langle a,b,| \ a=b^{-n}a^{-n}ba^{n}b^{n}, b=a^{-n}b^{-n}ab^{n}a^{n} 
\rangle \\
& = & \langle a,b,| \ ba^nb^n=a^nb^na, ab^na^n=b^na^nb  \rangle.
\end{eqnarray*}
So $G_n(T_r)$ and $G_n(T_l)$ are isomorphic for any $n$.\end{proof}

\section{\large \textbf{$G_n$ of the square and granny knots}}

The square knot $SK$ and the granny knot $GK$ are both connected sums of 
two trefoils; the granny knot is a connected sum of two right-handed 
trefoils or two left-handed trefoils, while the square knot is a connected 
sum of a right-handed trefoil and a left-handed trefoil.

It is well known that the fundamental groups of the complements of the
square knot and the granny knot are isomorphic; therefore, the knot group
and any knot invariants derived from it cannot distinguish these two knots. 
We will modify the defining presentations of $G_n(SK)$ and $G_n(GK)$ to a 
more symmetric form in this section. 

\begin{figure}[!ht]
\[
\includegraphics{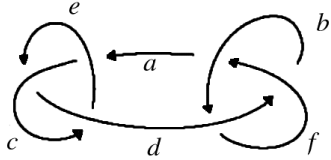}
\]
\caption{The square knot, SK}
\label{SK}
\end{figure}

From the diagram in Figure \ref{SK}, we obtain presentations of the 
fundamental quandle and $G_n$ of the square knot

\begin{eqnarray*}
Q(SK) & = & \langle a,b,c,d,e,f \ | \ a\tr e=c, \ e\tr c=d, \ c\tr d=e,\\
& & f\tr d=b, \ b\tr f=d, a\tr b=f \rangle,\\
G_n(SK) & = &  \langle a,b,c,d,e,f \ | \ a=e^nce^{-n}, \ e=c^ndc^{-n}, \\
& & 
\ c=d^ned^{-n},  f=d^nbd^{-n}, \ b=f^ndf^{-n}, a=b^nfb^{-n}. \rangle 
\end{eqnarray*}

\begin{proposition}
$G_n(SK)$ has presentation \[ \langle d,b,e \ | \ ed^ne^n=d^ne^nd,
\ bd^nb^n=d^nb^nd, \ e^nd^ned^{-n}e^{-n}=b^nd^nbd^{-n}b^{-n} \rangle.\]
\end{proposition}

\begin{proof}
We reduce $G_n(SK)$ to the required form via a sequence of Tietze moves. 
First, we eliminate the generator $a$:
\begin{eqnarray*}
G_n(SK) & = & \langle c,b,d,e,f \ | \ b^nfb^{-n}=e^nce^{-n}, \ e=c^ndc^{-n},
 c=d^ned^{-n}, \\ & &  f=d^nbd^{-n}, \ b=f^ndf^{-n}\rangle.
\end{eqnarray*}
Next, we eliminate the generator $c$:
\begin{eqnarray*}
G_n(SK) & = & \langle b,d,e,f \ | \ b^nfb^{-n}=e^n(d^ned^{-n})e^{-n}, \
e=(d^ne^nd^{-n})d(d^ne^{-n}d^{-n}), \\ & & f=d^nbd^{-n}, b=f^ndf^{-n} 
\rangle \\
& = & \langle e,d,b,f \ | \ b^nfb^{-n}=e^nd^ned^{-n}e^{-n},
\ e=d^ne^nde^{-n}d^{-n},  
\\ & & f=d^nbd^{-n}, \ b=f^ndf^{-n}\rangle. \end{eqnarray*}
Next, we eliminate the generator $f$:
\begin{eqnarray*}
G_n(SK) & = & \langle b,d,e \ | \ b^n(d^nbd^{-n})b^{-n}=e^nd^ned^{-n}e^{-n}, \
e=d^ne^nde^{-n}d^{-n},  \\
& & \ b=(d^nb^nd^{-n})d(d^nb^{-n}d^{-n}) \rangle\\
& = & \langle b,d,e \ | \ b^nd^nbd^{-n}b^{-n}=e^nd^ned^{-n}e^{-n}, \
e=d^ne^nde^{-n}d^{-n}, \\ & & b=d^nb^ndb^{-n}d^{-n}
\rangle.\end{eqnarray*}
Rewriting slightly, we obtain
\begin{eqnarray*}
G_n(SK) & = &\langle b,d,e \ | \ bd^nb^n=d^nb^nd, \ ed^ne^n=d^ne^nd, \\ & &
b^nd^nbd^{-n}b^{-n}=e^nd^ned^{-n}e^{-n} \rangle,\end{eqnarray*}
as required. \end{proof}

In the case $n=1,$ the relation \[b^nd^nbd^{-n}b^{-n}=e^nd^ned^{-n}e^{-n}\] 
becomes \[bdbd^{-1}b^{-1}=eded^{-1}e^{-1},\] which is 
a consequence of the other two relations, and
we obtain \[G_1(SK)= \pi_1(S^3\setminus SK)=\langle b,d,e \ | \ ded=ede,
\ dbd=bdb \rangle.\] Note that for $n>1$ we no longer have deficiency 1.

\begin{figure}[!ht]
\[
\includegraphics{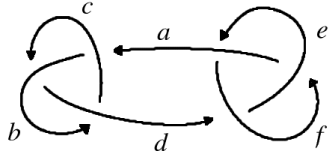}
\]
\caption{The granny knot, GK}
\label{GK}
\end{figure}

From the diagram in Figure \ref{GK}, we obtain presentations of the 
fundamental quandle and $G_n$ of the granny knot,
\begin{eqnarray*} Q(GK)  & = & \langle a,b,c,d,e,f \ | \  a\tr c=b, 
\ c\tr b=d, \ b\tr d=c, \ d\tr f=e, \\ & &  f\tr e=a, e\tr a=f \rangle \\
 G_n(GK) & = & \langle a,b,c,d,e,f \ | \ a=c^nbc^{-n}, \ c=b^ndb^{-n}, \
b=d^ncd^{-n}, \ d=f^nef^{-n}, \\ & & f=e^nae^{-n},  e=a^nfa^{-n} \rangle. 
\end{eqnarray*}

\begin{proposition}
$G_n(GK)$ has presentation 
\begin{eqnarray*}
G_n(GK)& = & 
\langle d, b, e \ | \ bd^nb^n=d^nb^nd,
\ de^nd^n=e^nd^ne, \\ & &  
\ e^{-n}d^{-n}ed^ne^n=b^nd^nbd^{-n}b^{-n}. \rangle 
\end{eqnarray*}
\end{proposition}

\begin{proof}
We reduce the presentation of $G_n(GK)$ obtained from the diagram. We begin by 
eliminating the generator $f$:

\begin{eqnarray*}
G_n(GK) & = & \langle a,b,c,d,e \ | \ a=c^nbc^{-n}, \ c=b^ndb^{-n}, \
b=d^ncd^{-n}, \\ 
& &  d=(e^na^ne^{-n})e(e^na^{-n}e^{-n}), \ e=a^n(e^nae^{-n})a^{-n} \rangle \\
& = &\langle a,b,c,d,e \ | \ a=c^nbc^{-n}, \ c=b^ndb^{-n}, \ b=d^ncd^{-n}, \\
& & d=e^na^nea^{-n}e^{-n},  ea^ne^n=a^ne^na \rangle \\
& = &\langle a,b,c,d,e \ | \ a=c^nbc^{-n}, \ c=b^ndb^{-n}, \ b=d^ncd^{-n}, \\
& & de^na^n=e^na^ne,  \ ea^ne^n=a^ne^na \rangle.\end{eqnarray*}

Next, eliminating the generator $c$ yields
\begin{eqnarray*}
G_n(GK) & = & \langle a,b,d,e \ | \ a=(b^nd^nb^{-n})b(b^nd^{-n}b^{-n}), 
\\ & & b=d^n(b^ndb^{-n})d^{-n},\ de^na^n=e^na^ne,  ea^ne^n=a^ne^na  \rangle \\
& = &\langle a,b,d,e \ | \ a=b^nd^nbd^{-n}b^{-n}, \ b=d^nb^ndb^{-n}d^{-n},
\\ & & de^na^n=e^na^ne, \ ea^ne^n=a^ne^na  \rangle \\
& = &\langle a,b,d,e| a=b^nd^nbd^{-n}b^{-n}, \ bd^nb^n=d^nb^nd,
\\ & & de^na^n=e^na^ne, \ ea^ne^n=a^ne^na  \rangle.\end{eqnarray*}

Eliminating the generator $a$, using the fact that
\[(b^nd^nbd^{-n}b^{-n})^m=b^nd^nb^md^{-n}b^{-n},\] yields
\begin{eqnarray*}
G_n(GK) & = & \langle b,d,e \ | \
bd^nb^n=d^nb^nd, \\ & & 
de^n(b^nd^nbd^{-n}b^{-n})^n=e^n(b^nd^nbd^{-n}b^{-n})^ne, \\
& & e(b^nd^nbd^{-n}b^{-n})^ne^n=(b^nd^nbd^{-n}b^{-n})^ne^n(b^nd^nbd^{-n}b^{-n})
\rangle \\
& = & \langle b,d,e \ | \ bd^nb^n=d^nb^nd, \\ & &
de^n(b^nd^nb^nd^{-n}b^{-n})=e^n(b^nd^nb^nd^{-n}b^{-n})e, \\ & &
e(b^nd^nb^nd^{-n}b^{-n})e^n=(b^nd^nb^nd^{-n}b^{-n})e^n(b^nd^nbd^{-n}b^{-n})
\rangle.
\end{eqnarray*}

The relation $bd^nb^n=d^nb^nd$ says that we can move $b$ past $d^nb^n$ to get
$d$. In particular, we can apply this relation $n$ times to obtain
$b^nd^nb^n=d^nb^nd^n$. Then
\begin{eqnarray*}
G_n(GK) & = & \langle b,d,e \ | \
bd^nb^n=d^nb^nd, \\ & & de^n(b^nd^nb^n)d^{-n}b^{-n}=e^n(b^nd^nb^n)d^{-n}b^{-n}e, \\ 
& & e(b^nd^nb^n)d^{-n}b^{-n}e^n=(b^nd^nb^n)d^{-n}b^{-n}e^nb^nd^nbd^{-n}b^{-n}
\rangle \\
& = & \langle b,d,e \ | \ 
bd^nb^n=d^nb^nd, \\ & & de^n(d^nb^nd^n)d^{-n}b^{-n}=e^n(d^nb^nd^n)d^{-n}b^{-n}e, \\ 
& & e(d^nb^nd^n)d^{-n}b^{-n}e^n=(d^nb^nd^n)d^{-n}b^{-n}e^nb^nd^nbd^{-n}b^{-n}
\rangle \\
& = & \langle b,d,e \ | \  bd^nb^n=d^nb^nd, de^nd^n=e^nd^ne, \\ & &
ed^ne^n=d^ne^nb^nd^nbd^{-n}b^{-n} \rangle \\
& = & \langle b,d,e \ | \  bd^nb^n=d^nb^nd, de^nd^n=e^nd^ne, \\ & & 
e^{-n}d^{-n}ed^ne^n=b^nd^nbd^{-n}b^{-n} \rangle,
\end{eqnarray*}
as required.\end{proof}

As with the square knot, in the case $n=1$ the relation 
$e^{-n}d^{-n}ed^ne^n=b^nd^nbd^{-n}b^{-n}$ is a consequence of the other two 
relations, and we obtain 
\[G_1(GK)=\langle b,d,e \ | \ bdb=dbd, ede=ded 
\rangle \cong G_1(SK).\] When $n>1$, we have 
\begin{eqnarray*}
G_n(SK) & = & \langle d, b, e \ | \ bd^nb^n=d^nb^nd, \ ed^ne^n=d^ne^nd, \\ 
& & e^nd^ned^{-n}e^{-n}=b^nd^nbd^{-n}b^{-n}. \rangle \end{eqnarray*}
and
\begin{eqnarray*}
G_n(GK) & = & \langle d, b, e \ | \ bd^nb^n=d^nb^nd, \ de^nd^n=e^nd^ne, \\ 
& & e^{-n}d^{-n}ed^ne^n=b^nd^nbd^{-n}b^{-n} \rangle. \end{eqnarray*}
 
\section{\large \textbf{Counting homomorphisms from $G_n(SK)$ and 
$G_n(GK)$ to a finite group}}

We will use the presentations of $G_n(SK)$ and $G_n(GK)$ 
given in the end of the last section. Denote
\[D=d^n,\ B=b^n, \ E=e^n.\]
Then the subgroup of $G_n(SK)$ generated by $D,E,B$ is isomorphic to 
\[G_1(SK)\cong\langle D,B,E\ | \ BDB=DBD,\ EDE=DED \rangle.\]
That this subgroup is isomorphic to $G_1(SK)$ follows from the topological 
interpretation of $G_n$ given in \cite{W}. We write the presentation of 
$G_n(SK)$ as follows:
\[G_n(SK)= \langle d, b, e \ | \ bDB=DBd, \ eDE=DEd, \  
EDeD^{-1}E^{-1}=BDbD^{-1}B^{-1} \rangle. \]
If we plug the first and the second relations into the third relation, the
later becomes
\[EDDEd(EDDE)^{-1}=BDDBd(BDDB)^{-1}.\]
Since $D=d^n$, we have $Dd=dD$ and the third relation is equivalent to
\[(ED)^3d(ED)^{-3}=(BD)^3d(BD)^{-3}.\]

Let $H$ be a finite group. Homomorphisms of $G_n(SK)$ into $H$ may be 
constructed in the following way. Suppose first that we are given a
homomorphism $\rho:G_1(SK)\longrightarrow H$. To define an extension 
$\hat{\rho}:G_n(SK)\longrightarrow H$, we first choose
$\hat\rho(d)\in H$ such that $\hat\rho(d)^n=\rho(D)$. Then define
\[\hat\rho(b)=\rho(D)\rho(B)\hat\rho(d)\rho(B)^{-1}\rho(D)^{-1}
\quad \mathrm{and} \quad
\hat\rho(e)=\rho(D)\rho(E)\hat\rho(d)\rho(E)^{-1}\rho(D)^{-1}.\]

Then if $\hat\rho$ also satisfies 
\[(\rho(E)\rho(D))^3\hat\rho(d)(\rho(E)\rho(D))^{-3}=
(\rho(B)\rho(D))^3\hat\rho(d)(\rho(B)\rho(D))^{-3}
\]
for every choice of $\hat\rho(d)\in H$ such that $\hat\rho(d)^n=\rho(D)$, 
say that $H$ has property $T(n,SK)$. If $H$ satisfies property $T(n,SK)$
then every homomorphism from the subgroup into $H$ extends to a homomorphism
from $G_n(SK)$ into $H$.

\begin{example}
\textup{If $H$ is abelian, then $H$ has property $T(n,SK)$ for every 
$n\in \mathbb{Z}$. If the exponent of $H$ is a multiple of $n$, 
then $D=B=E=1$, and $H$ has property $T(n,SK)$.}
\end{example}

Contrary to our original conclusion, not every group satisfies property
$T(n,SK)$; the map $\rho:G_1(SK)\to S_{24}$ given by
\begin{eqnarray*}
\rho(B) & = & (1,8,10,5,2,7,9,6)(15,17,24,19,16,18,23,20) \\ 
\rho(D) & = & (3,5,12,7,4,6,11,8)(15,17,24,19,16,18,23,20) \\ 
\rho(E) & = & (3,5,12,7,4,6,11,8)(13,20,22,17,14,19,21,18) \\
\end{eqnarray*}
with 
\[\hat{\rho}(d)=(3,15,5,17,12,24,7,19,4,16,6,18,11,23,8,20)\]
satisfies $\hat{\rho}(d)^2=\rho(D)$, 
$\rho(D)\rho(B)\rho(D)=\rho(B)\rho(D)\rho(B)$ and
$\rho(D)\rho(E)\rho(D)=\rho(E)\rho(D)\rho(E)$ but fails to
satisfy $(\rho(E)\rho(D))^3\hat\rho(d)(\rho(E)\rho(D))^{-3}=
(\rho(B)\rho(D))^3\hat\rho(d)(\rho(B)\rho(D))^{-3}.$ Hence
$S_{24}$ does not satisfy property $T(n,SK)$. \cite{T}

\begin{proposition} \label{1}
Let $H$ be a finite group which satisfies property $T(n,SK)$.
For $h\in H$, denote $\sqrt[n]{h}=\{g\in H\ | \ g^n=h\}$. Then we have 
\[\left|\text{\rm Hom}(G_n(SK),H)\right|=\sum_{\rho\in\text{\rm Hom}(G_1(SK),
H)}\left|\sqrt[n]{\rho(D)}\right|.\]
\end{proposition}

For the group $G_n(GK)$, we have a similar presentation
\[G_n(GK)= \langle d, b, e \ | \ bDB=DBd, \ dED=EDe, \  
E^{-1}D^{-1}eDE=BDbD^{-1}B^{-1} \rangle. \]
Thus, we can define an analogous property $T(n,GK)$, and proposition 
\ref{1} still holds if we replace $G_n(SK)$ by $G_n(GK)$. Therefore, 
we have the following conclusion.
 
\begin{corollary} \label{2}
For any finite group $H$ satisfying both properties $T(n,SK)$ and $T(n,Gk)$
we have 
\[
\left|\text{\rm Hom}(G_n(SK),H)\right|=\left|\text{\rm Hom}(G_n(GK),H)\right|
\]
for each $n$.
\end{corollary}

Note that $H$ acts on $\text{Hom}(G_n(SK),H)$ by conjugation. By our 
construction, two homomorphisms $\rho_1,\rho_2\in\text{Hom}(G_1(SK),H)$ are 
conjugate to each other iff their extensions 
$\hat\rho_1,\hat\rho_2\in\text{Hom}(G_n(SK),H)$ are conjugate to each other.
So we have similar results like proposition \ref{1} and corollary \ref{2} for 
the numbers 
of conjugacy classes of $\text{Hom}(G_n(SK),H)$ and $\text{Hom}(G_n(GK),H)$, 
respectively, for each $n$.   

\section{\large \textbf{Are $G_n(SK)$ and $G_n(GK)$ isomorphic for $n\geq 2$?}}

There is obviously a subtle difference in the presentations of $G_n(SK)$ and 
$G_n(GK)$ for $n\geq 2$. So even though we have not been able to distinguish 
these two groups via counting homomorphisms into finite groups, we think that 
these two groups should be not isomorphic to each other for each $n\ge 2$.

In \cite{W2}, Wada defined the twisted Alexander polynomial for a finitely 
presented group $G$ associated with two homomorphisms from $G$ into 
a free abelian group of finite rank and a matrix group, respectively. For the 
groups $G_n(SK)$ and $G_n(GK)$, their abelianizations are both infinite cyclic.
So we can have a twisted Alexander polynomial if we have a homomorphism into 
a matrix group. 
Using {\tt Maple}, we calculated the twisted Alexander polynomials
for $G_n(SK)$ and $G_n(GK)$ associated with homomorphisms into $SL(2,p)$ and 
$PSL(2,p)$ for some small primes $p$ and small $n$. So far, we are not able 
to distinguish these two groups using twisted Alexander polynomials. For 
example, using the isomorphism $PSL(2,7)\cong SL(3,2)$, we calculated the 
total of 8232 twisted Alexander polynomials for both $G_3(SK)$ and $G_3(GK)$ 
associated with homomorphisms into $PSL(2,7)$. It turns out that these two 
sets of 8232 twisted Alexander polynomials are equal to each other. 

We have also computed $|\mathrm{Hom}(G_n(SK),H)|$ and 
$|\mathrm{Hom}(G_n(GK),H)|$ for various values of $n$ for selected groups $H$
with orders as large as 360, always with the results that the two are equal.
Computing $|\mathrm{Hom}(G_n(SK),S_{24})|$ and $|\mathrm{Hom}(G_n(GK),S_{24})|$
is presently beyond the capabilities of our hardware.

Despite such unsuccessful attempts to show that $G_n(SK)$ and $G_n(GK)$ are 
not isomorphic for $n\ge 2$, we still propose the following conjecture:

\begin{conjecture} $G_n(SK)$ and $G_n(GK)$ are not isomorphic for each 
$n\geq 2$.
\end{conjecture}

\end{document}